\documentclass[10pt,a4paper,draft]{article}

\usepackage[leqno]{amsmath}
\usepackage{amsfonts,amssymb,amsthm}

\usepackage[all]{xy}
\SelectTips{cm}{12}   

\newcommand{\Db}{\mathrm{D^b}}
\newcommand{\DD}{{\mathrm{D}}}
\newcommand{\HH}{{\mathrm{H}}}
\newcommand{\LL}{{\mathbb{L}}}
\newcommand{\RR}{{\mathbb{R}}}
\newcommand{\FM}{\mathsf{FM}}
\newcommand{\MM}{\mathsf{M}}
\newcommand{\st}{\,:\,}
\DeclareMathOperator{\Aut}{Aut}
\DeclareMathOperator{\Spec}{Spec}
\DeclareMathOperator{\Hilb}{Hilb}
\newcommand{\HHilb}{{\rm{\text{-}Hilb}}}

\DeclareMathOperator{\Hom}{Hom}
\newcommand{\shHom}{\mathcal{H}om}
\DeclareMathOperator{\id}{id}
\DeclareMathOperator{\Sp}{Sp}
\newcommand{\can}{\mathsf{can}}      
\newcommand{\perm}{\mathsf{perm}}    
\newcommand{\abelian}{\mathrm{ab}}   
\newcommand{\asterisk}{\star}
\newcommand{\inv}{^{-1}}
\newcommand{\compose}{\asterisk}
\DeclareMathOperator{\Coh}{Coh}
\DeclareMathOperator{\Qcoh}{Qcoh}

\newcounter{enumcounter}
\newenvironment{enumliste}{\begin{list}{(\arabic{enumcounter})}
                      {\usecounter{enumcounter}
                       \setlength{\topsep}{0ex}
                       \setlength{\partopsep}{0ex}
                       \setlength{\listparindent}{0ex}
                       \setlength{\itemsep}{0ex}
                       \setlength{\parsep}{0ex}
                       \setlength{\leftmargin}{2em}
                       \setlength{\labelwidth}{1.5em}
                       \setlength{\parskip}{0ex}
                      }
                      }{\end{list}}

\newcommand{\isom}{ \text{{\hspace{0.48em}\raisebox{0.8ex}{${\scriptscriptstyle\sim}$}}}
                    \hspace{-0.65em}{\rightarrow}\hspace{0.3em}} 
\newcommand{\embed}{\hookrightarrow}
\newcommand{\surject}{\rightarrow\hspace{-1.8ex}\rightarrow}
\newcommand{\map}[1]{\stackrel{{}_#1}{\longrightarrow}}
\newcommand{\ratmap}{\dashrightarrow}
\newcommand{\sqmatnobraces}[4]{{\genfrac{.}{.}{0pt}{1}{#1}{#3} \:
                        \genfrac{.}{.}{0pt}{1}{#2}{#4}} }

\newcommand{\IC}{\mathbb C}
\newcommand{\IQ}{\mathbb Q}

\newcommand{\IZ}{\mathbb Z}
\newcommand{\ko}{\mathcal{O}}
\newcommand{\kt}{\mathcal{T}}
\newcommand{\kz}{\mathcal{Z}}

\newcommand{\bib}[5]{{\bibitem{#1} #2, {\emph{#3},} #4#5.}}
\newcommand{\Lin}{{\textrm{Lin}}_G}
\newcommand{\for}{{\mathit{for}}}           %
\renewcommand{\inf}{{\mathit{inf}}}         %
\newcommand{\foraut}{{\mathit{for}}}        %
\newcommand{\infaut}{{\mathit{inf}}}        %
\newcommand{\liftaut}{{\mathit{lift}}}      %
\newcommand{\inorlov}{{\eta}}               %
\newcommand{\sporlov}{{\gamma}}             %
\newcommand{\gdelta}{{G_\Delta}}            %
\newcommand{\Gex}{{\mathrm{ex}}}            %
\newcommand{\cohom}{{(\cdot)^\HH}}          %
\newcommand{\res}{{\mathit{res}}}
\newcommand{\BMeq}{\mathrm{eq}}             %

\newcounter{lemma}
\theoremstyle{plain}                        %
\newtheorem{localtheorem}[lemma]{Theorem}
\newtheorem{locallemma}[lemma]{Lemma}
\newtheorem{localproposition}[lemma]{Proposition}

\theoremstyle{definition}
\newtheorem{localremark}[lemma]{Remark}
\newtheorem{localremarks}[lemma]{Remarks}
\newtheorem{localexample}[lemma]{Example}

\begin{document}

\begin{center}
\textbf{\Large Equivariant autoequivalences for finite group actions}

\bigskip

David Ploog
\end{center}

\begin{quote}{\small\scshape Abstract}
The familiar Fourier-Mukai technique can be extended to an equivariant
setting where a finite group $G$ acts on a smooth projective variety $X$.
In this paper we compare the group of invariant autoequivalences
$\Aut(\Db(X))^G$ with the group of autoequivalences of $\DD^G(X)$.
We apply this method in three cases: Hilbert schemes on K3 surfaces,
Kummer surfaces and canonical quotients.
\end{quote}

\section{Introduction and Setup}

It often proves useful to consider analogues of classical settings,
adding the presence of a group action.
Instances of this in algebraic geometry are e.g.
equivariant intersection theory or the McKay type theorems. There already
is a theory for derived categories of varieties with actions by algebraic
groups \cite{bernstein-lunts}. In this article, we study the behaviour of
automorphism groups of such derived categories in the case when the group
is finite. This paper grew out of my Ph.D. thesis \cite{phd-ploog}. I 
would like to thank Georg Hein, Daniel Huybrechts, Manfred Lehn, and Richard 
Thomas for their help and valuable suggestions.

We always work with varieties over $\IC$. A kernel $P\in\Db(X\times Y)$
gives rise to a Fourier-Mukai transform which we denote by
$\FM_P:\Db(X)\to\Db(Y)$. We also introduce special notation for
composition of such transforms: Let us write
$\FM_Q\circ\FM_P=\FM_{Q\compose P}:\Db(X)\to\Db(Z)$ for $P\in\Db(X\times
Y)$ and $Q\in\Db(Y\times Z)$, i.e.
 $Q\compose P=\RR p_{XZ*}(p_{XY}^*P\otimes^\LL p_{YZ}^*Q)$.
This works for morphisms as well: $f:P\to P'$ and $g:Q\to Q'$ give rise 
to $g\compose f:Q\compose P\to Q'\compose P'$.

\subsection{Linearisations and $\DD^G(X)$}

Let $X$ be a smooth projective variety on which a finite group $G$ acts. A
\emph{$G$-linearisation} of a sheaf $E$ on $X$ is given by isomorphisms
$\lambda_g:E\isom g^*E$ for all $g\in G$ satisfying $\lambda_1=\id_E$ and
$\lambda_{gh}=h^*\lambda_g\circ\lambda_h$. A morphism
$f:(E_1,\lambda_1)\to(E_2,\lambda_2)$ is \emph{$G$-invariant}, if
 $f=g\cdot f:=\lambda_{2,g}\inv\circ g^*f \circ \lambda_{1,g}$ for all
$g\in G$.
The category of $G$-linearised coherent sheaves on $X$ with $G$-invariant
morphisms is denoted by $\Coh^G(X)$; note that it is abelian and contains
enough injectives, see \cite{bkr}. Put $\DD^G(X):=\Db(\Coh^G(X))$ for its
derived category.

There is an equivalent point of view on $\DD^G(X)$: let $\kt$ be the
category consisting of $G$-linearised objects of $\Db(X)$, i.e. complexes
$E^\bullet\in\Db(X)$ together with isomorphisms $\lambda_g:E^\bullet\isom
g^*E^\bullet$ in $\Db(X)$ satisfying the same cocycle condition as above.
The canonical functor $\DD^G(X)\to\kt$ is fully faithful in view of
 $\Hom_{\DD^G(X)}((E_1,\lambda_1)^\bullet,(E_2,\lambda_2)^\bullet) =
   \Hom_{\Db(X)}(E_1^\bullet,E_2^\bullet)^G$
for objects $(E_1,\lambda_1)^\bullet,(E_2,\lambda_2)^\bullet\in\DD^G(X)$.
To show that the functor is essentially surjective, take
$(E^\bullet,\lambda)\in\kt$. Choosing an injective bounded resolution
$E^\bullet\isom I^\bullet$, we can assume that $\lambda$ corresponds to
genuine complex maps $\tilde\lambda_g:I^\bullet\isom g^*I^\bullet$. Hence
$(I,\tilde\lambda)^\bullet$ is a complex with a linearisation of each
sheaf, and using
 $\Db(X) \cong \DD^{\mathrm{b}}_{\Coh(X)}(\Qcoh(X))
                \ni(I,\tilde\lambda)^\bullet$,
we get $\kt\cong\Db(X)$.

\subsection{Equivariant push-forwards and Fourier-Mukai transforms}

Now consider two such varieties with finite group actions, $(X,G)$ and
$(X',G')$. A map between them is given by a pair of morphisms
$\Phi:X\to X'$ and $\varphi:G\to G'$ such that
 $\Phi\circ g=\varphi(g)\circ\Phi$ for all $g\in G$. Then, we have the
pull-back $\Phi^*:\Coh^{G'}(X')\to\Coh^{G}(X)$ (and its derived functor
$\LL \Phi^*:\DD^{G'}(X')\to\DD^{G}(X)$) which just means
equipping the usual pull-back $\Phi^*E'$ with the $G_1$-linearisation
 $\Phi^*\lambda'_{\varphi(g)} :
   \phi^*E' \isom \Phi^*\varphi(g)^*E' = g^*\Phi^*E'$.

Suppose that $\varphi$ is surjective and put $K:=\ker(\varphi)$. Then,
there is also an equivariant push-forward defined for
$(E,\lambda)\in\Coh^G(X)$ in the following way: the usual push-forward
$\Phi_*E$ is canonically $G$-linearised since $\varphi$ is surjective. Now
$K$ acts trivially on $X'$, thus is it possible to take $K$-invariants.
Then the subsheaf $\Phi_*^K(E,\lambda) := [\Phi_*E]^K \subset \Phi_*E$
is still $G'$-linearised and $\RR\Phi_*^K:\DD^G(X)\to\DD^{G'}(X')$ is
the correct push-forward.

For objects $(E,\lambda)\in\Coh^G(X)$ and $(E',\lambda')\in\Coh^{G'}(X')$
and a $G$-invariant morphism $\Phi^*E'\to E$, the adjoint morphism
$E'\to\Phi_*E$ has image in $\Phi^K_*E$. Hence, the
functors $\Phi^*:\Coh^{G'}(X')\to\Coh^G(X)$ and
$\Phi_*^K:\Coh^G(X)\to\Coh^{G'}(X')$ are adjoint; analogously for
$\LL\Phi^*$ and $\RR\Phi^K_*$.

As a consequence, the
usual Fourier-Mukai calculus extends to the equivariant setting if we use
these functors (the tensor product of two linearised objects is obviously
again linearised). Explicitly, for an object
 $(P,\rho)\in\DD^{G\times G'}(X\times X')$
we get a functor
\[ \FM_{(P,\rho)} : \DD^G(X) \to \DD^{G'}(X'), \quad
   (E,\lambda) \mapsto \RR p_{X'}^G( P\otimes^\LL p_X^*(E)) \]
whit the projections $p_{X'}:X\times X'\to X'$ and $p_X:X\times X'\to X$
(and similar projections on the group level).

\subsection{Inflation and restriction}

There is an obvious forgetful functor
 $\for:\DD^G(X) \to \Db(X)$. In the other direction, we have the inflation
functor $\inf:\Db(X)\to\DD^G(X)$ with $\inf(E):=\bigoplus_{g\in G}g^*E$
and the $G$-linearisation comes from permuting the summands%
\footnote{
\ In a similar vein, every symmetric polynomial
$\sigma\in\IZ[x_1,\dots,x_n]^{S_n}$ (where $n:=\#G$)
gives rise to a functor $\inf_\sigma:\Db(X)\to\DD^G(X)$. For example,
$\inf_{x_1+\dots+x_n}=\inf$ and
 $\inf_{x_1\cdots x_n}(E)=\bigotimes_{g\in G}g^*E$.
As an application, for an ample line bundle $L$ on $X$ we get
$\inf_{x_1\cdots x_n}(L)$, an ample $G$-linearised line bundle.
}.
A generalisation of $\inf$ to the case of a subgroup $H\subset
G$ is given by
\[ \inf^G_H : \DD^H(X) \to \DD^G(X), \quad
              (E,\lambda) \mapsto \bigoplus_{[g]\in H\setminus G}g^*E \]
and the $G$-linearisation of the sum is a natural combination of $\lambda$
and permutations.

See Bernstein/Lunts \cite{bernstein-lunts} for generalisations of
$\DD^G(X)$ and $\inf^G_H$ to the case of algebraic groups (neither of
which is straight forward).

\subsection{Invariant vs linearised objects}

Obviously, a $G$-linearised object has to be $G$-invariant, i.e. fixed by
all pull-backs $g^*$. It is a difficult question under which conditions
the other direction is true. For us the following fact (\cite[Lemma
3.4]{phd-ploog}) will suffice.
\begin{locallemma} \label{linearisation-class}
Let $E\in\Db(X)$ be simple and $G$-invariant. Then there
is a group cohomology class $[E]\in\HH^2(G,\IC^*)$ such that $E$ is
$G$-linearisable if and only if $[E]=0$. Furthermore, if $[E]=0$, then
the set of $G$-linearisations of $E$ is canonically a $\hat{G}$-torsor.
\end{locallemma}

\begin{proof}
Note that the $G$-action on $\Aut(E)=\IC^*$ is trivial. There are
isomorphisms $\mu_g:E\isom g^*E$ for all $g\in G$. As $E$ is simple, we
can define units $c_{g,h}\in\IC^*$ by
 $\mu_{gh}=h^*\mu_g\circ\mu_h\cdot c_{g,h}$.
It is a straightforward check that the map $c:G^2\to\IC^*$ is a 2-cocycle
of $G$ with values in $\IC^*$, i.e. $c\in Z^2(G,\IC^*)$. Replacing the
isomorphisms $\mu_g$ with some other $\mu'_g$ yields the map $e:G\to\IC^*$
such that $\mu'_g=\mu_g\cdot e_g$. The two cocycles $c,c':G^2\to\IC^*$
derived from $\mu$ and $\mu'$ differ by the boundary coming from $e$ by
another easy computation. Hence, $c/c'=d(e)$ and thus
 $c=c'\in \HH^2(G,\IC^*)$.
Thus the $G$-invariant object $E$ gives rise to a unique class
$[E]:=c\in \HH^2(G,\IC^*)$. In these terms, $E$ is $G$-linearisable if and
only if $c\equiv1$, i.e.\ $[E]$ vanishes.

For the second statement, we write $\hat{G}:=\Hom(G,\IC^*)$ for the group
of characters and $\Lin(E)$ for the set of non-isomorphic
$G$-linearisations of $E$. Consider the $\hat{G}$-action
 $\hat{G}\times\Lin(E)\to\Lin(E)$,
 $(\chi,\lambda)\mapsto\chi\cdot\lambda$
on $\Lin(E)$. First take $\chi\in\hat{G}$ and $\lambda\in\Lin(E)$ such
that $\chi\cdot\lambda=\lambda$. Then, there is an isomorphism
$f:(E,\lambda)\isom(E,\chi\cdot\lambda)$ which in turn immediately implies
$\chi=1$ using $f\in\Aut(E)=\IC^*$. Thus, the action is effective.
Now take two elements $\lambda, \lambda'\in\Lin(E)$ and consider
$\lambda_g^{-1}\circ\lambda'_g:E\isom g^* E\isom E$. As $E$ is simple, we
have $\lambda_g^{-1}\circ\lambda'_g=\chi(g)\cdot\id_E$. It follows from
the cocycle condition for linearisations that $\chi$ is multiplicative,
i.e.\ $\chi\in\hat{G}$. In other words, $\lambda'=\chi\cdot\lambda$ and
the action is also transitive. Altogether $\hat{G}$ acts simply transitive
on $\Lin(E)$.
\end{proof}

For our use of group cohomology, refer for example \cite{serre}. The
second cohomology $\HH^2(G,\IC^*)$ of the finite group $G$ acting
trivially (on an algebraically closed field of characteristic 0) is also
known as the \emph{Schur multiplier} of $G$ (refer \cite[\S25]{huppert}).
Two relevant facts about it are:
$\HH^2(G,\IC^*)$ is a finite abelian group; its exponent is a divisor of
$\#G$. Examples are given by
$\HH^2((\IZ/n\IZ)^k,\IC^*)=(\IZ/n\IZ)^{k(k-1)/2}$
for copies of a cyclic group; $\HH^2(D_{2n},\IC^*)=\IZ/2\IZ$ and
$\HH^2(D_{2n+1},\IC^*)=0$ for the dihedral groups with $n>1$; and
$\HH^2(S_{n},\IC^*)=\IZ/2\IZ$ for the symmetric groups with $n>3$.

Note that a group with vanishing Schur multiplier has the following
property: every simple $G$-invariant object of $\DD^G(X)$ is
$G$-linearisable, no matter how $G$ acts on $X$.
\begin{localremark}
The condition that $E$ be simple in the Lemma is important.
Consider an abelian surface $A$ with the action of
$G=\IZ/2\IZ=\{\pm\id_A\}$. Then the sheaf $E:=k(a)\oplus k(-a)=\inf(k(a))$
is $G$-invariant but not simple. Yet it is uniquely
$\IZ/2\IZ$-linearisable as an easy computation shows
 \cite[Example 3.9]{phd-ploog}
(in contrast to $G$-invariant simple sheaves, which have precisely two
non-isomorphic $G$-linearisations according to $\HH^2(\IZ/2\IZ,\IC^*)=0$).
This behaviour is expected from geometry: by derived McKay correspondence
(\cite{bkr}) one has $\DD^G(A)\cong\Db(X)$, where $X$ is the Kummer
surface of $A$, a crepant resolution $\psi:X\to A/G$. Under this
equivalence, skyscraper sheaves of points $x\in X$ outside of exceptional
fibres of $\psi$ are mapped to $k(\psi(x))\oplus k(-\psi(x))$.
\end{localremark}

\begin{localexample} \label{canonical-linearisation}
If as before $G$ acts on $X$, then the canonical sheaf $\omega_X$ is
simple (as it is a line bundle) and $G$-invariant (because it is
functorial). Due to this functoriality, it is actually $G$-linearisable:
the morphism $g:X\to X$ induces a morphism of cotangent bundles
$g_*:g^*\Omega_X\to\Omega_X$. Going to determinants and using adjunction
yields the desired isomorphisms
 $\lambda_g:=\det(g\inv_*):\omega_X\isom g^*\omega_X$.
\end{localexample}

\section{Groups of autoequivalences}

We are interested in comparing the automorphism group $\Aut(\DD^G(X))$ with
the group
 $\Aut(\Db(X))^G :=
   \{F\in\Aut(\Db(X))\st g^*\circ F=F\circ g^* \ \forall g\in  G\}$.
It turns out that a useful intermediate step is to look at
Fourier-Mukai equivalences on $\Db(X)$ which are
diagonally $G$-linearised.

To make this precise, consider a Fourier-Mukai transform
$\FM_P:\Db(X)\to\Db(X')$. Suppose that $G$ acts on both $X$ and $X'$. Then
we have the diagonal action $G\times X\times X'\to X\times X'$,
$g\cdot(x,x'):=(gx,gx')$ which we sometimes (especially in the case
$X=X'$) for emphasis call the $G_\Delta$-action of $G$ on $X\times X'$.
Now we are in a position to study objects
 $(P,\rho)\in\DD^{G_\Delta}(X\times X')$ which give Fourier-Mukai
equivalences $\FM_P:\Db(X)\isom\Db(X')$. In other words, these are
ordinary kernels for equivalences $\Db(X)\isom\Db(X')$ which additionally
have been equipped with a $G_\Delta$-linearisation.

Not every kernel in $P\in\Db(X\times X')$ has the latter property. A
necessary condition is that $P$ must be $G_\Delta$-invariant, i.e.
$(g,g)^*P\cong P$ for all $g\in G$, or, equivalently,
$g^*\circ\FM_P=\FM_P\circ g^*$. Now we apply the following general fact:
\begin{locallemma} \label{simple-kernel}
If $P\in\Db(X\times Y)$ is the Fourier-Mukai kernel of an equivalence
$\FM_P:\Db(X)\isom\Db(Y)$ then $P$ is simple, i.e.
 $\Hom_{\Db(X\times Y)}(P,P)=\IC$.
\end{locallemma}
\begin{proof}
Fix $f\in\Hom_{\Db(X\times Y)}(P,P)$ and let $Q$ be a quasi-inverse
kernel for $\FM_P$, i.e. $P\compose Q\cong \ko_{\Delta_Y}$. By
 $\Hom_{\Db(Y\times Y)} (\ko_{\Delta_Y},\ko_{\Delta_Y}) =
   \Hom_{Y\times Y} (\ko_{\Delta_Y},\ko_{\Delta_Y}) = \IC$
we have $f\compose\id_Q=c\cdot\id_{\ko_{\Delta_Y}}$ for a
$c\in\IC$. Composing again with $\id_P:P\to P$ gives
$ f = f\compose\id_{\ko_{\Delta_X}} = f\compose (\id_Q\compose \id_P)
    = (f\compose \id_Q)\compose \id_P
    = (c\cdot\id_{\ko_{\Delta_Y}})\compose \id_P
    = c\cdot\id_P $.
\end{proof}

Combining Lemmas \ref{linearisation-class} and \ref{simple-kernel}, we
see that $G_\Delta$-invariant kernels for equivalences are
$G_\Delta$-linearisable, provided that $\HH^2(G,\IC^*)=0$ (or more
generally, if the obstruction class in $\HH^2(G,\IC^*)$ vanishes).

Now suppose we have an arbitrary object
 $(P,\rho)\in\DD^{G_\Delta}(X\times X')$
and the accompanying functor $\FM_P:\Db(X)\to\Db(X')$. The general device
of inflation allows us to produce the following equivariant Fourier-Mukai
transform from $(P,\rho)$:
\[ \FM_{(P,\rho)}^G := \FM_{\inf_{G_\Delta}^{G^2}(P,\rho)} :
    \DD^G(X)\to\DD^G(X') . \]
For brevity, we set
 $G\cdot P:=\inf_{G_\Delta}^{G^2}(P,\rho)\in\DD^{G^2}(X\times X')$
and call it \emph{the inflation of $(P,\rho)$}.
The following lemma states the main properties of this assignment.
\begin{locallemma} \label{inflated-kernel-properties}
Let $X$, $X'$, $X''$ be smooth projective varieties with $G$-actions and
let $(P,\rho)\in\DD^{G_\Delta}(X\times X')$ and 
$(P',\rho')\in\DD^{G_\Delta}(X'\times X'')$.
\begin{enumliste}
\item $\FM_{(\ko_{\Delta_X},\can)}^G\cong\id:\DD^G(X)\to\DD^G(X)$.
\item For $(P,\rho)\in\DD^{G_\Delta}(X\times X')$ there is a commutative
      diagram
      \[\xymatrix@C=4em{
      \DD^G(X) \ar[r]^{\FM^G_{(P,\rho)}} \ar[d]^{\for} & \DD^G(X') \ar[d]^\for \\
      \Db(X) \ar[r]^{\FM_P} & \Db(X')
       }\]
\item
      $\FM_{(P',\rho')}^G \circ \FM_{(P,\rho)}^G
                           \cong\FM_{(P'\compose P,\rho'\compose\rho)}^G$
      where $(\rho'\compose\rho)_g:=\rho'_g\compose\rho_g$.
\item
      $\FM_P$ fully faithful $\implies$ $\FM_{(P,\rho)}^G$ fully faithful.
\item
      $\FM_P$ equivalence $\implies$ $\FM_{(P,\rho)}^G$ equivalence.
\end{enumliste}
\end{locallemma}
\begin{proof}
(1) The structure sheaf $\ko_\Delta$ of the diagonal
 $\Delta\subset X\times X$
has a canonical $G_\Delta$-linearisation, as $(g,g)^*\ko_\Delta=\ko_\Delta$.
The inflation of $\ko_\Delta$ is
 $G\cdot\ko_\Delta = \bigoplus_{g\in G} \ko_{(g,1)\Delta}$,
and its $G^2$-linearisation is given by the permutation of
summands via $G\to G$, $g\mapsto kgh\inv$.

Using this one can check by hand that $\FM^G_{(\ko_\Delta,\can)}$ maps any
$(E,\lambda)\in\DD^G(X)$ to itself (see \cite[Example 3.14]{phd-ploog}).

(2) Take any object $(E,\lambda)\in \DD^G(Y)$. Then, we have by
definition of equivariant Fourier-Mukai transforms
 $\FM^G_P : \DD^G(Y) \to \DD^G(Y)$,
 $(E,\lambda) \mapsto [\RR p_{2*}(G\cdot P\otimes^\LL p_1^*E)]^{G\times 1}.$
The $G\times1$-linearisation of $G\cdot P$ is given by permutations (the
$G_\Delta$-linearisation of $P$ does not enter). Since
 $\RR p_{2*}((g,1)^*P\otimes^\LL p_1^*E) =
  \RR p_{2*}(g,1)^*(P\otimes^\LL p_1^*{g\inv}^*E) =
  \RR p_{2*}(P\otimes^\LL p_1^*E)$,
we see that
 $\RR p_{2*}(G\cdot P\otimes^\LL p_1^*E) \cong
   \bigoplus_G\RR p_{2*}(P\otimes^\LL p_1^*E)$
and $G\times1$ acts with permutation matrices where the $1$'s are replaced
by $p_1^*\lambda_g$'s. Taking $G\times1$-invariants singles out a
subobject of this sum isomorphic to one summand.

A morphism $f:E_1\to E_2$ in $\DD^G(X)$ is likewise first taken to a
$G$-fold direct sum. The final taking of $G\times1$-invariants then leaves
one copy of $\FM_P(f)$.

(3) The composite $\FM_{G\cdot P'}\circ\FM_{G\cdot P}$ has the kernel
\begin{align*}
  (G\cdot P')\compose (G\cdot P)
& =
  [\RR p_{13*}(p_{12}^*(G\cdot P)\otimes^\LL p_{23}^*(G\cdot P'))]^{1\times G\times 1} \\
& =
  [\RR p_{13*}(p_{12}^*\bigoplus_{g\in G}(g,1)^*P \otimes^\LL
           p_{23}^*\bigoplus_{h\in G}(h,1)^*P') ]^{1\times G\times 1} \\
& \cong_{\rho'}
  [\RR p_{13*}(p_{12}^*\bigoplus_{g\in G}(g,1)^*P \otimes^\LL
           p_{23}^*\bigoplus_{h\in G}(1,h\inv)^*P') ]^{1\times G\times 1}
           \\
& =
  [\bigoplus_{g,h\in G}\RR p_{13*} (g,1,h\inv)^*
          (p_{12}^*P \otimes^\LL p_{23}^*P') ]^{1\times G\times 1} \\
& =
  [\bigoplus_{g,h\in G} (g,h\inv)^* \RR p_{13*}
          (p_{12}^*P \otimes^\LL p_{23}^*P') ]^{1\times G\times 1} .
\end{align*}
Now note that
$(1,c,1)\in 1\times G\times1$ acts on $(G\cdot P')\compose(G\cdot P)$ via
permutations (inverse multiplications from left) and $\rho$ on $P$, and
$(1,c,1)$ acts purely by permutations (which are multiplications from right)
on $P'$.
Plugging this into the last equation, we find that after taking invariants
we end up with
 $\bigoplus_{d\in G}(d\inv,d)^*\RR p_{13*} (p_{12}^*P\otimes^\LL p_{23}^*P')$.
Since the $(d\inv,d)$'s give all classes in $G_\Delta\setminus G^2$, we
find that $(G\cdot P')\compose(G\cdot P)\cong G\cdot (P'\compose P)$.

(4) Fix two objects $(E_1,\lambda_1),(E_2,\lambda_2)\in \DD^G(X)$. The
injectivity of the natural map
 $\Hom_{\DD^G(X)}((E_1,\lambda),(E_2,\lambda_2)) \to
   \Hom_{\DD^G(X')}(\FM_P^G(E_1,\lambda_1),\FM_P^G(E_2,\lambda_2))$
is a consequence of two facts:
 $\Hom_{\DD^G(X)}(\cdot,\cdot) = \Hom_{\Db(X)}(\cdot,\cdot)^G
                        \subset \Hom_{\Db(X)}(\cdot,\cdot)$
by the definition of morphisms in $\DD^G(X)$ on one hand and 
 $\Hom_{\Db(X)}(E_1,E_2) \cong \Hom_{\Db(X')}(\FM_P(E_1),\FM_P(E_2))$
by hypothesis.

Similarly, the surjectivity uses the same facts. One just replaces the
embedding
 $\Hom_{\Db(X)}(\cdot,\cdot)^G \subset \Hom_{\Db(X)}(\cdot,\cdot)$
with the averaging projection (Reynolds operator)
 $\theta:\Hom_{\Db(X)}(E_1,E_2) \surject \Hom_{\Db(X)}(E_1,E_2)^G$
given by
 $\theta(f) :=
  \frac{1}{\#G}\sum_{g\in G} \lambda_{2,g}\inv\circ g^*f\circ\lambda_{1,g}$.

(5) follows from (3) and (1): let $\FM_P:\Db(X)\isom \Db(X')$ be an
equivalence. Then,
 $Q=\RR\shHom(P,\ko_{X\times X'})\otimes p_X^*\omega_X[\dim(X)]$
is the Fourier-Mukai kernel of a quasi-inverse for $\FM_P$. As the
canonical bundle $\omega_X$ is $G$-linearisable (see Example
\ref{canonical-linearisation}), $Q$ inherits a $G_\Delta$-linearisation
from those of $P$ and $p_X^*\omega_X$. We have
 $(G\cdot P)\compose(G\cdot Q)=G\cdot(P\compose Q)=G\cdot\ko_\Delta$.
There are precisely $\#\hat{G}$ different $G$-linearisations for
$\omega_X$ as well as $\#\hat{G}$ different $G_\Delta$-linearisations for
$\ko_\Delta$. Thus, exactly one choice of $G$-linearisation for $\omega_X$
will equip the composition $P\compose Q=\ko_\Delta$ with the canonical
$G_\Delta$-linearisation. But then (1) shows that this is the kernel of
the identity on $\DD^G(X')$. Hence, $G\cdot P$ is an equivalence kernel as
was $P$.
\end{proof}

Let us consider the following automorphism groups:
\begin{align*}
  \Aut(\Db(X))^{G}        & =  \{ \FM_P\in \Aut(\Db(X)) \st
                                  (g,g)^*P \cong P \: \: \forall g\in G \} \\
  \Aut(\DD^G(X))          & =  \{ \FM_{\widetilde{P}}:\DD^G(X)\isom \DD^G(X) \st
                                  \widetilde{P}\in\DD^{G^2}(X\times X) \} \\
  \Aut^{G_\Delta}(\Db(X)) & :=  \{ (P,\rho)\in\DD^{G_\Delta}(X\times X) \st
                                  \FM_P\in\Aut(\Db(X)) \} .
\end{align*}
The first identity uses the action $G\times\Aut(\Db(X))\to\Aut(\Db(X))$
given by $g\cdot F:= (g\inv)^*\circ F\circ g^*$ and the formula
$\FM_{(g,g)^*P}=g^*\circ\FM_P\circ (g\inv)^*$ and finally Orlov's result
on the existence of Fourier-Mukai kernels \cite[Theorem 2.2]{orlov-k3}.
This has been extended by Kawamata to smooth stacks associated to normal
projective varieties with quotient singularities \cite{kawamata-stacks}.
In view of $\Coh([X/G])\cong\Coh^G(X)$, this implies the second relation.
Finally we have to turn $\Aut^{G_\Delta}(\Db(X))$ into a group. This is
done by Lemma \ref{inflated-kernel-properties}: (3) settles the
composition, (1) the neutral element and (5) the inverses.

The following theorem \cite[Proposition 3.17]{phd-ploog} is an attempt
to compare these groups.
\begin{localtheorem} \label{AutDG}
Let the finite group $G$ act on a smooth projective variety $X$.\
\begin{enumliste}
\item The construction of inflation gives a group homomorphism $\infaut$
      which fits in the following exact sequence, where $Z(G)\subset G$ is
      the centre of $G$:
\[\xymatrix@1@R=1ex{
 0 \ar[r] & Z(G) \ar[r] & \Aut^{G_\Delta}(\Db(X)) \ar[r]^-\infaut & \Aut(\DD^G(X)) \\
          &             & ~(P,\rho)~ \ar@{|->}[r]                & ~\FM_{(P,\rho)}^G~
      } \]
\item Forgetting the $G_\Delta$-linearisation gives a group homomorphism
      $\foraut$ which sits in the following exact sequence; here
      $G_\abelian:=G/[G,G]\cong\Hom(G,\IC)^*=\HH^1(G,\IC^*)$ is the 
      abelianisation:
      \begin{equation*} \xymatrix@1@R=1ex{
       0 \ar[r] & G_\abelian \ar[r] &
       \Aut^{G_\Delta}(\Db(X)) \ar[r]^-\foraut & \Aut(\Db(X))^G \ar[r] & \HH^2(G,\IC^*) \\
       & & ~(P,\rho)~ \ar@{|->}[r] & ~\FM_P~ \ar@{|->}[r] & ~[P]~
      } \end{equation*}
\end{enumliste}
\end{localtheorem}

\begin{proof}
(1) It follows from Lemma \ref{inflated-kernel-properties} that
$\infaut$ is a group homomorphism. The kernel $\ker(\infaut)$ consists
of $(P,\rho)\in\DD^{G_\Delta}(X^2)$ giving equivalences such that
$G\cdot P\cong G\cdot\ko_\Delta$. Obviously, this forces $P$ to be a
sheaf of type $P\cong(g,1)^*\ko_\Delta$ for some $g\in G$. Now
$(g,1)^*\ko_\Delta$ is $G_\Delta$-invariant if and only if $g\in Z(G)$ as
 $ (h,h)^*(g,1)^*\ko_\Delta = (gh,h)^*\ko_\Delta \cong
   (gh,h)^*(h\inv,h\inv)^*\ko_\Delta = (h\inv gh,1)^*\ko_\Delta $.
This in turn implies that the isomorphism $(g,1):X^2\isom X^2$ is a
$G_\Delta$-map. In particular, $P\cong(g,1)^*\ko_\Delta$ gets the pulled
back $G_\Delta$-linearisation. Giving $\ko_\Delta$ a
$G_\Delta$-linearisation $\lambda\in\hat{G}$ different from the canonical
one yields
 $G\cdot(\ko_\Delta,\lambda)\not\cong G\cdot(\ko_\Delta,\can)$;
this follows for example from the uniqueness of Fourier-Mukai kernels.
Both facts together imply $\ker(\infaut)\cong Z(G)$.

(2) It is obvious from the definition of
$\Aut^{G_\Delta}(\Db(X))$ that $\foraut$ is a group homomorphism. The
kernel of $\foraut$ corresponds to the $G_\Delta$-linearisations on
$\ko_\Delta$. From Lemma \ref{linearisation-class}, we see that they form
a group isomorphic to $\hat{G}=\Hom(G,\IC^*)\cong G_\abelian$. Note also
$\Hom(G,\IC^*)=\HH^1(G,\IC^*)$ as $G$ acts trivial on $\IC^*$. Given
$\FM_P\in\Aut(\Db(X))^{G_\Delta}$, we know from Lemma \ref{simple-kernel}
that its Fourier-Mukai kernel $P$ is simple. Furthermore, it is
$G_\Delta$-invariant by assumption, so that the map $\FM_P\mapsto[P]$ is
defined as in Lemma \ref{linearisation-class}. To see that it is a group
homomorphism, take two $G_\Delta$-invariant kernels $P,Q\in\Db(X^2)$.
Choose isomorphisms $\lambda_g:P\isom(g,g)^*P$ and $\mu_g:Q\isom(g,g)^*Q$
for all $g\in G$. Then we have
$\lambda_{gh}=(h,h)^*\lambda_g\circ\lambda_h\cdot[P]_{g,h}$ and likewise
for $Q$, $\mu$. Furthermore, the composition of $\lambda_g$ and $\mu_g$
gives an isomorphism
 $\mu_g\compose\lambda_g : Q\compose P \isom ((g,g)^*Q)\compose((g,g)^*P)$
and the latter term is canonically isomorphic to $(g,g)^*(Q\compose P)$.
Then
 $\mu_{g,h}\compose\lambda_{gh} =
  (h,h)^*(\mu_g\compose\lambda_g) \circ (\mu_h\compose\lambda_h)
                                 \cdot [Q\compose P]_{g,h}$.
The formula
 $(B\circ A)\compose(D\circ C)\cong(B\compose D)\circ(A\compose C)$
for morphisms $P''\map{B}P'\map{A}P$ and $Q''\map{D}Q'\map{C}Q$ implies
$[Q\compose P]_{g,h}=[Q]_{g,h}\cdot[P]_{g,h}$. Now it is obvious that
$[\cdot]\circ\foraut=0$ and finally $\FM_P\in\Aut(\Db(X))^G$ with $[P]=0$
implies that $P$ is $G_\Delta$-linearisable by Lemma
\ref{linearisation-class}.
\end{proof}

\section{Applications}

The applications will make use of the derived McKay correspondence: assume
that $X$ is a smooth (quasi-)projective variety with an action of the
finite group $G$. Then there is a $G$-Hilbert scheme $G\HHilb(X)$ which
parametrises $G$-clusters in $X$, i.e. points $\xi\in G\HHilb(X)$
correspond to $0$-dimensional subschemes $Z_\xi\subset X$ such that
$\HH^0(\ko_{Z_\xi})\cong\IC[G]$, as $G$-representations. Typical examples
of such clusters are free $G$-orbits. The formal definition of
$G\HHilb(X)$ uses that it represents the relevant functor. Let
$\widetilde{X}\subset G\HHilb(X)$ be the connected component which contains
the free orbits. Then there is a birational morphism (the Hilbert-Chow
map) $\widetilde{X}\to X/G$. Combined with the projection $X\to X/G$, this
yields a universal subscheme $\kz\subset\widetilde{X}\times X$; note that
canonically $\ko_\kz\in\DD^{1\times G}(\widetilde{X}\times X)$. We refer to
\cite{bkr} for
\begin{localtheorem}[Bridgeland-King-Reid] \label{bkr}
Suppose that $\omega_X$ is locally trivial in $\Coh^G(X)$ and
$\dim(\widetilde{X}\times_{X/G}\widetilde{X})\leq\dim(X)+1$. Then
$\FM_{\ko_\kz}:\Db(\widetilde{X})\isom\DD^G(X)$ is an equivalence.
\end{localtheorem}

\subsection{Hilbert schemes on K3 surfaces}

At first consider two smooth projective varieties $X$ and $Y$ and the
$n$-fold products $X^n$, $Y^n$ with their natural $S_n$-actions. Let
$P\in\Db(X\times Y)$ be the kernel of a Fourier-Mukai transform
$\FM_P:\Db(X)\to\Db(Y)$. Then, the exterior tensor product
 $P^{\boxtimes n}=p_1^*P\otimes\cdots\otimes p_n^*P\in\Db(X^n\times Y^n)$
yields the functor
 $F^n=\FM_{P^{\boxtimes n}}:\Db(X^n)\to\Db(Y^n)$.
Furthermore, $P^{\boxtimes n}$ has an obvious $(S_n)_\Delta$-linearisation
via permutation of tensor factors. Hence, using inflation, we get the
new functor
 $F^{[n]}:=\FM_{P^{\boxtimes n}}^{S_n}:\DD^{S_n}(X^n)\to\DD^{S_n}(Y^n)$.
If we restrict to the case $X=Y$ and autoequivalences
$F:\Db(X)\isom\Db(X)$, we get a group homomorphism by Theorem
\ref{AutDG}
\[ \Aut(\Db(X)) \to \Aut^{(S_n)_\Delta}(\Db(X^n)) \to \Aut(\DD^{S_n}(X)),
    \quad \FM_P \mapsto \FM_{(P^{\boxtimes n},\perm)}^{S_n} . \]

From now on we need the provision $\dim(X)=\dim(Y)=2$. It is
well-known that for surfaces $\Hilb_n(X)$ is a crepant resolution
of $X^n/S_n$. Furthermore, a theorem of Haiman states
 $\Hilb_n(X)\cong S_n\HHilb(X^n)$, see \cite{haiman}.
If we additionally assume $\omega_X\cong\ko_X$ and $\omega_Y\cong\ko_Y$,
then we can invoke Theorem \ref{bkr} in order to obtain
 $\Db(\Hilb_n(X))\cong\DD^{S_n}(X^n)$
because in this case $X^n$ and $Y^n$ are symplectic manifolds and the
condition $\dim(\Hilb_n(X)\times_{X^n/S_n}\Hilb_n(X))<1+n\dim(X)$ of 
Theorem \ref{bkr} is automatically fulfilled; see 
\cite[Corollary 1.3]{bkr}. However, a posteriori this inequality is 
true for general surfaces as the dimension of the fibre product is a 
local quantity which may be computed with any (e.g. affine or 
symplectic) model. The above 
homomorphism of groups of autoequivalences is now
\[ \Aut(\Db(X)) \to \Aut(\Db(\Hilb_n(X))) . \]
It is always injective: this is clear for $n>2$ since the centre of $S_n$
is trivial in this case. For $n=2$ one can check that the sheaf
$\ko_{\Delta_X\times\Delta_X}$ with the non-canonical
$(S_2)_\Delta$-linearisation is not in the image of
$\Aut(\Db(X))\to\Aut^{(S_2)_\Delta}(\Db(X^2))$.

Let us introduce the shorthand
$X^{[n]}:=\Hilb_n(X)$, so that $\Db(X^{[n]})\cong\DD^{S_n}(X^n)$ by the
above and $\Aut(\Db(X))\embed\Aut(\Db(X^{[n]}))$.
The above technique shows
\begin{localproposition}
If $X$ and $Y$ are two smooth projective surfaces with $\Db(X)\cong\Db(Y)$,
then $\Db(X^{[n]})\cong\Db(Y^{[n]})$. \qed
\end{localproposition}

\begin{localremark} \label{hilbert-D-functorial}
A birational isomorphism $f:X\ratmap Y$ of smooth projective 
surfaces induces a birational map
$f^{[n]}:X^{[n]}\ratmap Y^{[n]}$ between their Hilbert schemes.
There is a derived analogue: a Fourier-Mukai transform (resp.\ equivalence) 
$F=\FM_P:\Db(X){\to}\Db(Y)$ induces a functor (resp.\ equivalence)
 $F^{[n]}=\FM_{P^{\boxtimes n}}^{S_n}:\Db(X^{[n]})\to\Db(Y^{[n]})$.
Since for the time being it is unknown whether every functor 
$F:\Db(X)\to\Db(Y)$ is of Fourier-Mukai type, we have to restrict to 
those (which include equivalences by \cite[Theorem 2.2]{orlov-k3}). 
\end{localremark}

Finally, we specialise to K3 surfaces.
\begin{localproposition} \label{birational-hilberts}
Let $X_1$ and $X_2$ be two projective K3 surfaces. If there is a
birational isomorphism $X_1^{[n]}\ratmap X_2^{[n]}$ of their Hilbert
schemes, then the derived categories are equivalent:
$\Db(X_1^{[n]})\cong\Db(X_2^{[n]})$.
\end{localproposition}

\begin{proof}
A birational isomorphism $f:X_1^{[n]}\ratmap X_2^{[n]}$ induces an
isomorphism on second cohomology,
 $f^*:\HH^2(X_1^{[n]})\isom \HH^2(X_2^{[n]})$, respecting the Hodge
structures, because Hilbert schemes of symplectic surfaces are
symplectic manifolds. From the crepant resolution
$X_1^{[n]}\to X_1^n/S_n$, we find
$\HH^2(X_1)\subset \HH^2(X_1^{[n]})$, and only the
exceptional divisor $E_1\subset X_1^{[n]}$ is missing:
$\HH^2(X_1^{[n]})=\HH^2(X_1)\oplus\IZ\cdot\delta_1$ with $2\delta_1=E_1$.
In particular, as $[E_1]$ is obviously an algebraic class, the
transcendental sublattices coincide:
$T(X_1)=T(X_1^{[n]})$. Hence, the birational isomorphism furnishes an
isometry $T(X_1)\cong T(X_2)$. Orlov's derived Torelli theorem for K3
surfaces \cite{orlov-k3} then implies $\Db(X_1)\cong \Db(X_2)$. But now
we apply the above result on lifting equivalences from K3 surfaces to
Hilbert schemes and deduce $\Db(X_1^{[n]})\cong\Db(X_2^{[n]})$, as
claimed.
\end{proof}

\begin{localremarks} (1) follows from Remark
\ref{hilbert-D-functorial} and (2) from Proposition
\ref{birational-hilberts}.\
\begin{enumliste}
\item
$\DD$-equivalent abelian or K3 surfaces $X_1$ and $X_2$ have
$\DD$-equivalent Hilbert schemes $X_1^{[n]}$ and $X_2^{[n]}$.
\item
Considering only birational equivalence classes of Hilbert schemes on
K3 surfaces, we find that each such class is finite, since K3 surfaces
have only finitely many Fourier-Mukai partners \cite{bm1}.
\item
Markman \cite{markman-brill-noether} gives an example of non-birational
Hilbert schemes $X_1^{[n]}$ and $X_2^{[n]}$ with
$\HH^2(X_1^{[n]})\cong\HH^2(X_2^{[n]})$. The above arguments still yield
$\Db(X_1)\cong\Db(X_2)$ and $\Db(X_1^{[n]})\cong\Db(X_2^{[n]})$, i.e.
$X_1^{[n]}$ is $\DD$-equivalent to $X_2^{[n]}$.
\end{enumliste}
\end{localremarks}

\subsection{Kummer surfaces}

Let $A$ be an abelian variety and consider the action of
 $G:=\{\pm\id_A\}\cong\IZ/2\IZ$.
In order to investigate $\Aut(\DD^G(A))$, we start with the exact
sequence (see Orlov's article \cite{orlov-abel}):
\[ \xymatrix{
  0 \ar[r] & \IZ\times A\times\hat{A} \ar[r]^\inorlov &
     \Aut(\Db(A)) \ar[r]^\sporlov &
      \Sp(A\times\hat{A}) \ar[r] & 0 }, \]
the first morphism $\inorlov$ maps a triple $(n,a,\xi)$ to the
autoequivalence $t_a^*\circ\MM_\xi[n]$, where $t_a:A\isom A$ denotes
the translation by $a$ and $\MM_\xi:\Db(A)\isom\Db(A)$ the line bundle
twist with $\xi$. Note that shifts, translations and twists by degree
0 line bundles commute.
Before turning to the second morphism $\sporlov$, we set
\[ \Sp(A\times\hat{A}) := \left\{
   \Bigl(\sqmatnobraces{f_1}{f_2}{f_3}{f_4}\Bigr)
           \in \Aut(A\times\hat{A}) ~:~
   \Bigl(\sqmatnobraces{f_1}{f_2}{f_3}{f_4}\Bigr)
   \Bigl(\sqmatnobraces{\phantom{-}\hat{f}_4}{-\hat{f}_2}{-\hat{f}_3}{\phantom{-}\hat{f}_1}\Bigr)
 = \Bigl(\sqmatnobraces{\id_A}{0}{0}{\id_{\hat{A}}}\Bigr) \right\} .
\]
Now given $F\in\Aut(\Db(X))$, there is a functorial way to attach
an equivalence $\Phi_F:\Db(A\times\hat{A})\isom\Db(A\times\hat{A})$
which sends skyscraper sheaves to skyscraper sheaves. Hence
$\Phi_F$ yields an autmorphism
 $\gamma(F):A\times\hat{A}\isom A\times\hat{A}$
which turns out to be in $\Sp(A\times\hat{A})$
(see the original \cite[\S2]{orlov-abel} by Orlov or
\cite[\S4]{phd-ploog} for a slightly different presentation).

Note that $G$ acts on $\Aut(\Db(A))$ via conjugation, i.e.
 $(-1)\cdot F:=(-1)^*\circ F\circ(-1)^*$.
This induces an action on $\Sp(A\times\hat{A})$, which is trivial since
$\sporlov((-1)^*)=-\id_{A\times\hat{A}}$ is central. Taking $G$-invariants
of Orlov's exact sequence, we get
\[ \xymatrix{
  0 \ar[r] & \IZ\times A[2]\times\hat{A}[2] \ar[r]^-{\inorlov^G} &
     \Aut(\Db(A))^G \ar[r]^-{\sporlov^G} &
      \Sp(A\times\hat{A}) \ar[r] & 0 }. \]
Here, $A[2]\subset A$ and $\hat{A}[2]\subset\hat{A}$ denote the
2-torsion subgroups; the surjectivity of $\sporlov^G$ uses
 $\HH^1(G,\IZ\times A\times\hat{A})=0$; see
  \cite[Proposition 4.8]{phd-ploog}.
Hence, any autoequivalence $F\in\Aut(\Db(A))$ differs from a
$G$-invariant one just by translations and degree 0 line bundle
twists.

Assume from now on that $A$ is an abelian \emph{surface}. Let $X$ be the
corresponding Kummer surface, which is a crepant resolution of $A/G$. A
realisation is $X=G\HHilb(A)$ and we use derived McKay correspondence
(Theorem \ref{bkr}) to infer $\Db(X)\cong\DD^G(A)$. From Theorem
\ref{AutDG} we get group homomorphisms
\[ \xymatrix@1@C=3em{
 \Aut(\Db(A))^G &
  \Aut^\gdelta(\Db(A)) \ar@{->>}[l]_-\foraut^-{2:1} \ar[r]^-\infaut_-{2:1}
 & \Aut(\DD^G(A))
} \]
In our situation both $\foraut$ and $\infaut$ are 2:1, and $\foraut$
is surjective as $\HH^2(\IZ/2\IZ,\IC^*)=0$. There does not seem to be a
homomorphism between the lower groups making the triangle commutative.

However, when going to cohomology (here with $\IQ$ coefficients
throughout) the roof can be completed to a diagram (refer
\cite[Proposition 4.14]{phd-ploog}):
\[ \xymatrix{
       & \Aut(\Db(A)) \ar[r]^-\cohom & \Aut(\HH^{2*}(A)) \\
  \Aut^\gdelta(\Db(A)) \ar[ur]^{\foraut} \ar[dr]_{\infaut} \\
       & \Aut(\Db(X)) \ar[r]_-{\cohom} & \Aut(\HH^{*}(X)) \ar[uu]_\res
} \]
Here, $(\cdot)^\HH:\Aut(\Db(A))\to\Aut(\HH^{2*}(A))$ sends a Fourier-Mukai
equivalence to the corresponding isomorphism on cohomology. Further we use
that the image of $(\cdot)^\HH\circ\infaut$ lies inside the subgroup of
isometries preserving the exceptional classes,
\[ \Aut(\HH_\Gex^*(X)) :=
    \{ \varphi\in\Aut(\HH^*(X)) \st
        \varphi(\Lambda)=\Lambda \}
\]
where $\Lambda\cong\IQ^{16}$ is the lattice spanned by the $(-2)$-classes
arising from the Kummer construction; the morphism
$\res:\varphi\mapsto\varphi|_\Lambda$ is then the obvious restriction.

\subsection{Canonical quotients}

Let $X$ be a smooth projective variety whose canonical bundle is of finite
order. Suppose that $n>0$ is minimal with $\omega_X^n\cong\ko_X$. Then
there is an \'etale covering $\widetilde{X}\map{\pi}X$ of degree $n$ with
$\omega_{\widetilde{X}}\cong\ko_{\widetilde{X}}$. One concrete definition is
$\widetilde{X}=\Spec(\ko_X\oplus\omega_X\cdots\oplus\omega_X^{n-1})$, see
\cite[\S I.17]{bpvdv}. The group $G:=\IZ/n\IZ$ then acts freely on
$\widetilde{X}$ with $\widetilde{X}/G=X$. Fix a generator $g\in G$ and note that a
$G$-linearisation for some $\widetilde{E}\in\Db(\widetilde{X})$ is completely
determined by an isomorphism $\lambda:\widetilde{E}\isom g^*\widetilde{E}$ subject
to
 $(g^{n-1})^*\lambda\circ\cdots\circ g^*\lambda\circ\lambda=\id_{\widetilde{E}}$
as $G$ is cyclic. Here we have an equivalence
$\Coh(X)\cong\Coh^G(\widetilde{X})$ already on the level of abelian
categories (see \cite[\S7]{git-1st}). Hence,
$\Db(X)\cong\DD^G(\widetilde{X})$ as well, a fact which also follows from
the derived McKay correspondence using the trivially crepant
resolution $\xymatrix@1{X\ar[r]_{\id_X}&X}$. Then,
\[ \xymatrix@1@C=3em{
 \Aut(\Db(\widetilde{X}))^G &
  \Aut^\gdelta(\Db(\widetilde{X})) \ar@{->>}[l]_-\foraut^-{n:1} \ar[r]_-{n:1}^-\infaut
 & \Aut(\Db(X))=\Aut(\DD^G(\widetilde{X})) .
} \]
Bridgeland and Maciocia consider in \cite{bm3} canonical quotients from
the point of view of derived categories. They introduce the set of all
\emph{equivariant} equivalences by
\[ \Aut_\BMeq(\Db(\widetilde{X})) :=
    \{ (F,\mu) \in\Aut(\Db(\widetilde{X}))\times\Aut(G)
     \st g_*\circ\widetilde{F}\cong\widetilde{F}\circ\mu(g)_*
      \:\: \forall g\in G \};
\]
this is actually a group. There is an exact sequence
\[ \xymatrix{ 0 \ar[r] & \Aut(\Db(\widetilde{X}))^G \ar[r]
            & \Aut_\BMeq(\Db(\widetilde{X})) \ar[r]
            & \Aut(G) \ar[r] & 0
} \]
where the latter morphism maps $(\widetilde{F},\mu)\mapsto\mu$ and
$\Aut(\Db(\widetilde{X}))^G$ is by definition the group of all equivariant
equivalences with $\mu=\id_G$. Note that $G\cong\IZ/n\IZ$ implies
$\Aut(G)\cong\IZ/\varphi(n)\IZ$. Furthermore, we have a subgroup
$G\embed\Aut(\Db(\widetilde{X}))^G$, $g\mapsto g_*=(g\inv)^*$, or also
$G\embed\Aut_\BMeq(\Db(\widetilde{X}))$, $g\mapsto(g_*,\id_G)$. The latter 
is a normal divisor in view of
 $   (\widetilde{F},\mu)\inv\circ(g_*,\id_G)\circ(\widetilde{F},\mu)
   = (\widetilde{F}\inv,\mu\inv)\circ(g_*\circ\widetilde{F},\mu)
   = (\mu(g)_*,\id_G)$.
By \cite[Theorem 4.5]{bm3} every equivalence $F\in\Aut(\Db(X))$
has an equivariant lift $\widetilde{F}\in\Aut(\Db(\widetilde{X}))$, i.e.
 $\pi_*\circ\widetilde{F} \cong F\circ\pi_*$ and
 $\pi^*\circ F \cong \widetilde{F}\circ\pi^*$.
If $\widetilde{F}_1$ and $\widetilde{F}_2$ are two lifts of $F$, then
$\widetilde{F}_2\inv\circ\widetilde{F}_1$ is a lift of $\id_{\Db(X)}$ and thus
$\widetilde{F}_2\inv\circ\widetilde{F}_1\cong g_*$ for a $g\in G$ (\cite[Lemma
4.3(a)]{bm3}). Thus the equivariant lift
$\widetilde{F}\in\Aut_\BMeq(\Db(\widetilde{X}))$ is unique up to the action of $G$
and we get a group homomorphism
$\liftaut:\Aut(\Db(X))\to\Aut_\BMeq(\Db(\widetilde{X}))/G$. \cite[Lemma
4.3(b)]{bm3} states that if $F,F'\in\Aut(\Db(X))$ both lift to
$\widetilde{F}\in\Aut_\BMeq(\Db(\widetilde{X}))$, then they differ by a line
bundle twist: $F\cong F'\circ\MM_{\omega_X^i}$ with $0\leq i<n$. Thus
$\liftaut$ is $n:1$, and we propose the following commutative pentagram
\[ \xymatrix@!0@C=6em@R=6em{
   & \Aut^\gdelta(\Db(\widetilde{X}))
      \ar@{->>}[dl]_{\foraut}^{n:1} \ar[dr]_{n:1}^{\infaut} \\
       \Aut(\Db(\widetilde{X}))^G \ar@{^{(}->}[d]^{1:\varphi(n)} & &
        \Aut(\Db(X))=\Aut(\DD^G(\widetilde{X})) \ar[d]^\liftaut_{n:1} \\
   \Aut_\BMeq(\Db(\widetilde{X})) \ar@{->>}[rr]^{n:1} & &
   \Aut_\BMeq(\Db(\widetilde{X}))/G
} \]
The commutativity of this diagram boils down to the following
question: given a kernel $(P,\rho)\in\Aut^\gdelta(\Db(\widetilde{X}))$, is
$\FM_P:\Db(\widetilde{X})\isom\Db(\widetilde{X})$ a lift of
$\FM_{(P,\rho)}^G:\DD^G(\widetilde{X})\isom\DD^G(\widetilde{X})$ (where we identify
$\DD^G(\widetilde{X})\cong\Db(X)$)? However, this is clear from
$\pi_*:\Db(\widetilde{X})\to\DD^G(\widetilde{X})$, $E\mapsto\inf(E)$ and
$\pi^*:\DD^G(\widetilde{X})\to\Db(\widetilde{X})$, $(F,\lambda)\mapsto[F,\lambda]^G$.

\begin{localexample}
In the case of a double covering (e.g.\ $X$ an Enriques surface), we have
$n=2$ and hence $\Aut(\Db(\widetilde{X}))^G = \Aut_\BMeq(\Db(\widetilde{X}))$.
Then the diagram looks like
\[ \xymatrix@R=5em{
  & \Aut^\gdelta(\Db(\widetilde{X})) \ar[ld] \ar[dr]^\infaut \\
     \Aut(\Db(\widetilde{X}))^G/G
  & & \Aut(\Db(X)) \ar[ll]^\liftaut
} \]
\end{localexample}

\end{document}